\documentclass[12pt,a4paper,draft]{article}
\usepackage{amssymb,amsmath,stmaryrd,vmargin}
\usepackage[all]{xy}
\newtheorem{theo}{Theorem}[section]
\newtheorem{prop}[theo]{Proposition}

\newtheorem{lem}[theo]{Lemma}

\def \demdu#1 { {\sl Proof #1.} }

\def \Gal {\text{\rm Gal}}

\def \tor {\text {\rm Tor}}
\def \coker {\Coker}
\def \dem{{\sl Proof.} }

\def \limpro{\lim\limits_{\leftarrow} }

\def \qed{\text {$\quad \square$}}

\def \C {\overline {C}}
\def \U {\overline {U}}

\def \F#1 {\overline {F^\times_{#1}}}

\def\BM{{\mathbb{B}}}

\def\NM{{\mathbb{N}}}

\def\QM{{\mathbb{Q}}}

\def\ZM{{\mathbb{Z}}}

\def\XG{{\mathfrak X}}

\def\al{\alpha}

\def\ga{\gamma}
\def\Ga{\Gamma}

\def\la{\lambda}
\def\La{\Lambda}

\def\si{\sigma}

\def\th{\theta}

\def\ze{\zeta}

\def\DC{{\mathcal{D}}}

\def\GC{{\mathcal{G}}}
\def\HC{{\mathcal{H}}}

\def\NC{{\mathcal{N}}}

\def\UC{{\mathcal{U}}}

\def\UCt{{\widetilde{\mathcal{U}}}}

\def\Cti{{\widetilde{C}}}

\def\Mti{{\widetilde{M}}}

\def\Uti{{\widetilde{U}}}

\def\Aba{{\overline{A}}}

\def\Cba{{\overline{C}}}

\def\Iba{{\overline{I}}}

\def\Mba{{\overline{M}}}

\def\Uba{{\overline{U}}}

\def\Coker{\mathop{\mathrm{Coker}}\nolimits}

\def\Gal{\mathop{\mathrm{Gal}}\nolimits}

\def\im{\mathop{\mathrm{Im}}\nolimits}

\def\Ker{\mathop{\mathrm{Ker}}\nolimits}
\def\ker{\Ker}

\def\pr{\prime}

\def\MI{{\rm Minv}}

\providecommand{\MR}{\relax\ifhmode\unskip\space\fi MR }

\providecommand{\href}[2]{#2}
\def\Dbar{\leavevmode\lower.6ex\hbox to 0pt{\hskip-.23ex
    \accent"16\hss}D}
\def\cfac#1{\ifmmode\setbox7\hbox{$\accent"5E#1$}\else
    \setbox7\hbox{\accent"5E#1}\penalty 10000\relax\fi\raise 1\ht7
    \hbox{\lower1.15ex\hbox to 1\wd7{\hss\accent"13\hss}}\penalty 10000
    \hskip-1\wd7\penalty 10000\box7}
\def\cftil#1{\ifmmode\setbox7\hbox{$\accent"5E#1$}\else
    \setbox7\hbox{\accent"5E#1}\penalty 10000\relax\fi\raise 1\ht7
    \hbox{\lower1.15ex\hbox to 1\wd7{\hss\accent"7E\hss}}\penalty 10000
    \hskip-1\wd7\penalty 10000\box7}

\def\Dbar{\leavevmode\lower.6ex\hbox to 0pt{\hskip-.23ex \accent"16\hss}D}
  \def\cfac#1{\ifmmode\setbox7\hbox{$\accent"5E#1$}\else
  \setbox7\hbox{\accent"5E#1}\penalty 10000\relax\fi\raise 1\ht7
  \hbox{\lower1.15ex\hbox to 1\wd7{\hss\accent"13\hss}}\penalty 10000
  \hskip-1\wd7\penalty 10000\box7}
  \def\cftil#1{\ifmmode\setbox7\hbox{$\accent"5E#1$}\else
  \setbox7\hbox{\accent"5E#1}\penalty 10000\relax\fi\raise 1\ht7
  \hbox{\lower1.15ex\hbox to 1\wd7{\hss\accent"7E\hss}}\penalty 10000
  \hskip-1\wd7\penalty 10000\box7}

\def\ken {\Ker_n}
\def\Ken {\ken}
\def\Char {{\rm Char}}
\title{ Global Units modulo Circular Units :
descent without Iwasawa's Main Conjecture.}
\author{Jean-Robert Belliard\\
\ \\
\small{\it Universit\'e de Franche-Comt\'e,
Laboratoire de math\'ematiques UMR 6623,} \\
\small{\it 16 route de Gray,
25030 Besan\c con cedex,
FRANCE.}\\ \ \\
\small{\texttt{belliard@math.univ-fcomte.fr}}}
\begin{document}
\date{\today}
\maketitle
\setcounter{section}{-1}
\begin{abstract}  Iwasawa's classical asymptotical formula relates the orders of the $p$-parts $X_n$ of the ideal
class groups along a $\ZM_p$-extension $F_\infty/F$ of a number
field $F$, to Iwasawa structural invariants $\la$ and $\mu$ attached to the inverse limit $X_\infty=\limpro X_n$. It relies on "good" descent properties satisfied by
$X_n$.  If $F$ is abelian and $F_\infty$ is cyclotomic it is known
that the $p$-parts of the orders of the global units modulo
circular units $U_n/C_n$ are asymptotically equivalent to the
$p$-parts of the ideal class numbers. This suggests that these
quotients $U_n/C_n$, so to speak unit class groups, satisfy also
good descent properties. We show this directly, i.e. without using Iwasawa's Main Conjecture.
\end{abstract}

\section{Introduction}

Let $K$ be a number field and $p$ an odd prime ($p\neq 2$) and let
$K_\infty/K$ be a $\ZM_p$-extension (quite soon $K_\infty/K$ will
be the cyclotomic $\ZM_p$-extension). Recall the usual notations~:
$\Ga=\Gal(K_\infty/K)$ is the Galois group of $K_\infty/K$, $K_n$
is the $n^{\text{th}}$-layer  of $K_\infty$ (so that
$[K_n:K]=p^n$), $\Ga_n=\Gal(K_\infty/K_n)$, and
$G_n=\Gal(K_n/K)\cong \Ga/\Ga_n$. Let us consider a sequence
$(M_n)_{n\in \NM}$ of $\ZM_p[G_n]$-modules equipped with norm maps
$M_n\longrightarrow M_{n-1}$ and the inverse limit of this
sequence $M_\infty=\limpro M_n$ seen as a
$\La=\ZM_p[[\Ga]]$-module. The general philosophy of Iwasawa
theory is to study the simpler $\La$-structure of $M_\infty$, then
to try and recollect information about the $M_n$'s themselves from
that structure. For instance, if $M_\infty$ is $\La$-torsion, one
can attach two invariants $\la$ and $\mu$ to $M_\infty$. If we
assume further that
\begin{itemize}
\item[$(i)$] the $\Ga_n$ coinvariants
$(M_\infty)_{\Ga_n}$ are finite,
\item[$(ii)$] the sequence $(M_n)_{n\in
\NM}$ behaves "nicely" viz descent;
\end{itemize}
then one can prove that the
orders of the $M_n$ are asymptotically equivalent to $p^{\la n+\mu
p^n}$. These two assumptions are expected to occur whenever one
chooses for $(M_n)$ significant modules (from the number theoretic
point of view). However proving them may require some effort :
actually, $(i)$ for the canonical Iwasawa module
$\rm{tor}_\La(\XG_\infty)$ is (one of the many equivalent
formulations of) Leopoldt's conjecture. The historical example
occurs when we specialize $M_n=X_n$, the $p$-part of ideal class
group of $K_n$. Then the asymptotic formulas are a theorem of
Iwasawa. The proof of this theorem uses an auxiliary module
$Y_\infty$ which is pseudo-isomorphic to $X_\infty$ but with
better descent properties.

In the present paper we are interested in similar statements for
unit class groups, that is for $M_n=\Uba_n/\Cba_n$, the $p$-part of
the quotient of units modulo the subgroup of circular units of
$K_n$. In order that the latter merely exist we need to assume that all $K_n$ are abelian over $\QM$~:
hence $K$ is abelian and $K_\infty/K$ is cyclotomic. Let us assume further that
$K$ is totally real, which is not a loss of generality as long as
we are only interested in $p\neq2$. Now by Sinnott index formulas
we know that the orders of $\Uba_n/\Cba_n$ are asymptotically
equivalent to the orders of $X_n$. And as a consequence of
Iwasawa's Main Conjecture, the structural invariants $\la$ and
$\mu$ of both inverse limits $\Uba_\infty/\Cba_\infty=\limpro \U_n/\C_n$ and of $X_\infty=\limpro X_n$ are
equal. Using these two theorems we have immediately a somewhat indirect proof of an
analogue of Iwasawa theorem for unit class groups. Clearly a
direct proof of this fact must exist and the first goal of the present paper is to write it down.
The theorem \ref{iwaclassunit} shows that the orders of the unit class groups $\U_n/\C_n$ along the finite steps of the $\ZM_p$-extension are those prescribed by the structural Iwasawa  invariants of the inverse limit  $\U_\infty/\C_\infty$. It makes no use of any precise link between ideal and unit class groups. E.g. it does not use the Main Conjecture, neither even Sinnott's index formula.
On the other hand if one does use some classical results together with the theorem \ref{iwaclassunit}, one gets easier proofs of  beautiful and well known theorems (see \S \ref{applications}). For instance Sinnott's index formulas together with theorem \ref{iwaclassunit} imply the equality of the $\la$ and $\mu$ invariants of the two class groups {\it without using Iwasawa Main Conjecture}. The equality between these two $\la$ invariants is  actually  an important step (some times called "class number trick") in the proof of Iwasawa's Main Conjecture. If we use further Ferrero-Washington's theorem, we prove that all the $\mu$-invariants involved here are trivial.
Maybe some ideas in the present approach could be used in a
different framework where the equality of the two characteristic
ideals (that is Iwasawa's Main Conjecture) is still open.
In the paper \cite{Ng05}, \S 5, T. Nguyen Quang Do has also proven that
unit classes have asymptotically good descent properties, by
using Iwasawa's Main Conjecture in its strongest form : that is, following
the path outlined just above.

We conclude this introduction by recalling the (now) traditional
notations of cyclotomic Iwasawa theory. For any number field $F$
we put $S=S(F)$ for the set of places of $F$ dividing $p$. We will
adopt the following notations :
\begin{itemize}
\item[] For any abelian group $A$, we'll denote by $\Aba$ the $p$-completion of $A$, i.e.
the inverse limit  $\Aba=\underset \leftarrow \lim\ A/A^{p^m}$. If $A$ is finitely
generated over $\ZM$, then $\bar A \cong A \otimes \ZM_p$.
\item[] $U_F$ is the group of units of $F$.
\item[] $U^\prime_F$  is the group of $S$-units of $F$ : that is, elements of $F^\times$ whose valuations
are trivial for all finite places $v$ of $F$ such that $v\not\in S$.
\item[] $X_F$ is the $p$-part of $cl(F)$ which in turn is the ideal class group of $F$.
\item[] $X^\pr_F$ is the $p$-part of  $cl^\pr(F)$ which in turn is
the quotient of $cl(F)$ modulo the subgroup generated by classes of primes in $S$.
\item[] $\XG_F$ is the Galois group over $F$ of the maximal
 $S$-ramified (i.e. unramified outside $S$) abelian pro-$p$-extension of
$F$.
\item[] $\NC_F$ is the multiplicative group of the semi-local numbers. As a mere $\ZM_p$-module
$\NC_F= \prod_{v\in S} \overline{F_v^\times}$ where $F_v$ is the completion of $F$
at the place $v\in S$.
\item[] $\UC_F$ is the group of semi-local units of $F$. As a mere $\ZM_p$-module
$\UC_F=\prod_{v\in S} U^1_{v}$
where $U^1_{v}$ is the set of principal units of $F_v$, that is units $\equiv 1$ modulo the maximal ideal of $F_v$.

\end{itemize}
 So for instance the $\Uba_n$'s above could have been understood
as $U_{K_n}\otimes \ZM_p$
and $X_F=\overline {cl(F)}$.

\section{Consequences of Leopoldt's Conjecture}
If $F$ is Galois over $\QM$, the groups $\UC_F$ and $\NC_F$ come
equipped with the induced action of $\Gal(F/\QM)$. In other words
if we fix one place $v$ in $S$ and if we put $G_v=\Gal(F_v/\QM_p)$
then we have (as $\ZM_p$-modules)~: $$\UC_F\cong  U^1_{v}
\otimes_{\ZM_p[G_v]} \ZM_p[\Gal(F/\QM)]\ \text{and}\ \NC_F \cong
\overline {F_v^\times} \otimes_{\ZM_p[G_v]} \ZM_p[\Gal(F/\QM)].$$
These isomorphisms define the action of $\Gal(F/\QM)$ on $\UC_F$
and on $\NC_F$. Let us now consider the cyclotomic
$\ZM_p$-extension $K_\infty/K$ of our totally real abelian over
$\QM$ number field $K$ and the $n^{\text{th}}$-finite layers $K_n/K$,
i.e. $K_n$ is the unique subfield of $K_\infty$ with
$[K_n:K]=p^n$. We will denote by $C_n$ (and of course we'll be
more interested in $\Cba_n$) the group of circular units of $K_n$
as defined by Sinnott (\cite{Si80}). We will indicate consistently
by the subscript $_n$ arithmeticall object related to $K_n$. So
$\Uba_n$, $\UC_n$ and so on make sense. For any extension $L/F$ of
global fields the formula
$$N_{L/F} ((x_w)_{w\mid p})
=( \prod_{w\mid v} N_{L_w/K_v} (x_w))_{v\mid p}$$ defines a Galois
equivariant morphism $N_{L/F} \colon \NC_L\longrightarrow \NC_F$
which is compatible with the usual norms on global units for
instance. As ever $\Uba_\infty$, $\Cba_\infty$ $\NC_\infty$,
$\UC_\infty$, $\Uba^\pr_\infty$ and so on denotes the inverse
limit (related to norm maps) of $\Uba_n$, $\Cba_n$, $\NC_n$,
$\UC_n$, $\Uba^\pr_n$ and so on. Before considering more precise
properties of descent kernels (and cokernels) of unit classes, we
need to first establish their finiteness. This can surely be
extracted from various and older references as a consequence of
Leopoldt's conjecture (which is true here since all $K_n$ are
abelian over $\QM$). However a more precise and general result on
quotients of semi-local units modulo circular units can be found in
\cite{T99}. The point for our present approach is that the proofs
of \cite{T99} don't make any use of the Main Conjecture, and only
need Coleman morphisms.

We will need and freely use the following consequence of theorem 3.1 in \cite{T99}.
\begin{theo}\label{Leo}
Let $G=\Gal(K/\QM)$. Recall that $\Ga=\Gal(K_\infty/K)$ and
 that $\La=\ZM_p[[\Ga]]$. Fix a generator $\ga$ of $\Ga$.
For all non-trivial Dirichlet characters of the first kind $\psi$ of $\widehat G$,
following \cite{T99} let us denote by $g_\psi(T)$ the Iwasawa
power series associated to the Kubota-Leopoldt $p$-adic
$L$-function $L_p(\psi, s)$ (with our fixed choice of $\ga$).
\begin{enumerate}
\item Assume $p$ is tamely ramified in $K$ so that every characters of $G$ are of the first kind.
Then up to a power of $p$, the characteristic ideal of the $\La$-torsion module
$\UC_\infty/\Cba_\infty$ is generated by the product
$\displaystyle \prod_{\psi\in\widehat G,\psi\neq 1} g_\psi(T)$.
\item For all $n\in \NM$, the $\Ga_n$-coinvariants
$ (\UC_\infty/\Cba_\infty)_{\Ga_n}$ and $(\Uba_\infty/\Cba_\infty)_{\Ga_n}$ are finite.
\end{enumerate}
\end{theo}
\dem Let us consider $\psi$-parts $M^\psi$ as defined in
\cite{T99} for any $\ZM_p[G]$-module $M$. These $\psi$ parts are
naturally $\ZM_p[\psi]$-modules and $\ZM_p[G]$-submodules of $M$.
For technicalities about $\psi$-components  which are often left to readers see Beliaeva's thesis \cite{Tan}. From theorem 3.1 of
\cite{T99} we use that $\UC_\infty^\psi/\Cba_\infty^\psi$ is
$\ZM_p$-torsion free and that its characteristic series over
$\La[\psi]$ is $g_\psi(T)$. Let $\si\in\Gal(\QM_p[\psi]/\QM_p)$;
then $(\UC_\infty/\Cba_\infty)^{\psi^\si}$ is isomorphic to
$(\UC_\infty/\Cba_\infty)^\psi$ and
$g_{\psi^\si}(T)=\si(g_\psi(T))$. For $R$ equal  to $\La$ or to $\La[\psi]$, let us
abbreviate by $\Char_R(M)$ the characteristic ideal of the
$R$-module $M$. It is an easy linear algebra exercise to check
that if a power series $f$ generates $\Char_{\La[\psi]} (M)$ then
$N_{\QM_p[\psi]/\QM_p} (f)$ generates $\Char_\La(M)$. Hence we get
$$\Char_\La((\UC^\psi_\infty/\Cba_\infty^\psi)) =(N_{\QM_p[\psi]/\QM_p} (g_\psi(T)))
= \left (\prod_{\si \in \Gal(\QM_p[\psi]/\QM_p)} g_{\psi^\si}
(T)\right ).$$ Let us fix a set $\Psi$ of representatives of
$\widehat G$ up to $\QM_p$-conjugation classes. If a
$\ZM_p[G]$-module $M$ is $\ZM_p$-torsion free then the submodules
of $M$, $(M^\psi)_{\psi\in\Psi}$, are mutually direct summands and
the quotient $M/\oplus_{\psi\in\Psi} M^\psi$ is annihilated by
$\#G$. This applies to $\UC_\infty$ and $\Cba_\infty$, so that we
may use the inclusions $\bigoplus_{\psi\in\Psi} \UC_\infty^\psi
\subset \UC_\infty$ and $\bigoplus_{\psi\in\Psi}
\Cba^\psi_\infty\subset \Cba_\infty$ to form the snake diagram :
$$\xymatrix{0\ar[r] &\bigoplus_{\psi\in\Psi} \Cba_\infty^\psi \ar[r] \ar@{^(->}[d] &
 \bigoplus_{\psi\in\Psi} \UC^\psi_\infty \ar[r] \ar@{^(->}[d] & \bigoplus_{\psi\in\Psi}
 \UC_\infty^\psi/\Cba_\infty^\psi\ar[r]\ar@{-->}[d]^-\al & 0\\
 0\ar[r] & \Cba_\infty \ar[r] \ar[d]&\UC_\infty \ar[r]\ar[d] & \UC_\infty/\Cba_\infty\ar[r]\ar[d] &0 \\
 & \ZM_p-\text{torsion} \ar[r] & \ZM_p-\text{torsion}\ar[r]  &\coker\ \al \ar[r] &0 }$$
\noindent Now the kernel of $\al$ must be $\ZM_p$-torsion, as a
submodule of the cokernel of the first inclusion, hence is trivial
because all $\UC_\infty^\psi/\Cba_\infty^\psi$ are $\ZM_p$-free.
This proves that up to a power of $p$ the characteristic ideal
$\Char_\La(\UC_\infty/\Cba_\infty)$ is equal to
$\Char_\La(\bigoplus_{\psi\in\Psi} \UC_\infty/\Cba_\infty)$, which
in turn is what we wanted.

The second claim is equivalent to the fact that both characteristic ideals are prime to all $(T+1)^{p^n}-1$ for all $n\in\NM$. From the inclusion
$\Uba_\infty/\Cba_\infty \subset \UC_\infty/\Cba_\infty$ we see
that the first characteristic ideal divides the second.
Let $\ze$ be a $p$-power order root of unity and let $\rho$ be the unique
character of the second kind such that $\rho(\ga)=\ze^{-1}$. Then we have
$g_{\psi}(\ze-1)=L_p(1,\psi\rho)$. By Leopoldt's conjecture (see \cite{Wa}
corollary 5.30)  $L_p(1,\psi\rho)\neq 0$.
Hence in the special case where $p$ is tame in $K$, 2 follows from 1.
We now want to remove this hypothesis. Let $I_p\subset \Gal(K_\infty/\QM)$ be the
inertia subgroup of $p$, and let $\Iba_p\subset I_p$ be its
pro-$p$-part. Then the field $L=K_\infty^{\Iba_p}$ is a finite
abelian over $\QM$ number field and $p$ is at most tamely ramified
in $L$. Since all modules
$M_\infty$ depend on $K_\infty$ and not on $K$ itself, to prove
that the hypothesis is almost no loss of generality we just have
to show the following lemma~:
\begin{lem}\label{modram}
$L_\infty=K_\infty$.
\end{lem}
\dem Let $\BM_\infty/\QM$ be the $\ZM_p$-extension of $\QM$.
By construction we have $L\subset K_\infty$ and therefore
$L_\infty\subset K_\infty$. By definition of $L$, the group
$\Gal(K_\infty/L_\infty)$ is then a subgroup of $\Iba_p$, hence
$K_\infty/L_\infty$ is a $p$-extension totally ramified at $p$.
But on the other hand $\BM_\infty\subset L_\infty$ and $K_\infty$
is  abelian over $\QM$. Let $\QM^{\rm{ab}}$ be the maximal abelian
extension of $\QM$. Then $\QM^{\rm{ab}}/\BM_\infty$ is tamely
ramified at (the only) prime of $\BM_\infty$ above $p$. This shows
$L_\infty=K_\infty$.

\qed

Now, with that lemma we already have proven the finiteness
assertion when we take coinvariants
 along $\widetilde {\Ga}_n=\Gal(K_\infty/L_n)$. We conclude the proof of the second part of theorem
 \ref{Leo} by pointing that even if
the first few $\Ga_n$ may differ from $\tilde{\Ga}_n$, we get
$\Ga_n=\tilde{\Ga}_n$ as soon as $K_\infty/K_n$ is totally ramified at $p$.
Due to canonical surjections $(M_\infty)_{\Ga_n}\twoheadrightarrow (M_\infty)_{\Ga_{n-1}}$,
this does not change anything to the finiteness of $(\UC_\infty/\Cba_\infty)_{\Ga_n}$.

\qed

\begin{itemize}
\item[] {\it Remark :}
As long as one is only concerned with characteristic ideals up to
a power of $p$ and finiteness of coinvariants, theorem \ref{Leo}
was well known long before as part of the folklore, and can be
tracked back to Iwasawa (\cite{I64}; see also \cite{Grei92}) in some cases.
\end{itemize}

Now that the finiteness of $(\UC_\infty/\Cba_\infty)_{\Ga_n}$
(resp. $(\Uba_\infty/\Cba_\infty)_{\Ga_n}$) has been proved, we
want to relate these modules with their  counterparts
$\UC_n/\Cba_n$ (resp. $\Uba_n/\Cba_n$) at finite levels. On the way
we will study descent for various other multiplicative Galois
modules.

\section{Background: easy part of descent}
Let $M_\infty=\limpro (M_n)$ be a $\La$-module.  Of course the
projections $M_\infty\longrightarrow M_n$ factor through
$(M_\infty)_{\Ga_n} \longrightarrow M_n$, but in general these
maps have non trivial kernels and cokernels.  Let us denote $\ken
(M_\infty)$ for the kernel of  $(M_\infty)_{\Ga_n}\longrightarrow
M_n$, $\Mti_n\subset M_n$ for its image, and
$\Coker_n(M_\infty)=M_n/\Mti_n$ for its cokernel. Let $l$ be a
prime number ($l=p$ is allowed) and $E$ be a finite extension of
$\QM_l$. Put $E_\infty/E$ for a $\ZM_p$-extension of $E$, and
$E_n$ for the $n^{\rm th}$ finite layers. Let us abbreviate
$L_E=\overline {E^\times}$, consistently
$L_n=\overline{(E_n)^\times}$, and $L_\infty=\limpro L_n$. For all
$n\in \NM\cup{\infty}$, local class field theory identifies $L_n$
with $\Gal(M_n/E_n)$, where $M_n$ is the maximal abelian
$p$-extension of $E_n$.
\begin{lem}\label{locdes}
The natural map $(L_\infty)_{\Ga_n} \longrightarrow L_n$  fits
into an exact sequence $$\xymatrix{0\ar[r]& (L_\infty)_{\Ga_n}
\ar[r]& L_n\ar[r]^-{\rm Artin}& \Gal(E_\infty /E_n)\ar[r]& 0.}$$
\end{lem}
\dem This is well known to experts. We follow the cohomological
short cut of \cite{Ng84}. The group $\Ga_n\simeq \ZM_p$ is pro-$p$-free hence we have $H^2(\Ga_n,\QM_p/\ZM_p)=0$. Consider the inflation-restriction
sequence associated to the extension of groups
$$\xymatrix{1\ar[r]& \HC \ar[r] & \GC_n\ar[r] & \Ga_n \ar[r] &
1,}$$ where $\GC_n$ is the absolute Galois group of $E_n$ :
$$\xymatrix{0\ar[r]& H^1(\Ga_n,\QM_p/\ZM_p)  \ar[r] &H^1(\GC_n,\QM_p/\ZM_p) \ar[r] & H^1(\HC,\QM_p/\ZM_p)^{\Ga_n}  \ar[r]
& 0.}$$  Applying Pontryagin duality
and class field theory we get
$$\xymatrix{0\ar[r] & (L_\infty)_{\Ga_n} \ar[r]&
L_n\ar[r]^-{\rm Artin}& \Gal(E_\infty /E_n)\ar[r]& 0.}$$

\qed

Let us come back to our global field $K$ and its cyclotomic $\ZM_p$-extension.
\begin{prop}\label{slocdes} Let $S_n$ be the set of places of $K_n$ dividing $p$. For all $n$ and $v\in S_n$
put $K_{n,v}$ for the completion at $v$ of $K_n$ and
$(K_{n,v})^{\rm cyc}_\infty/K_{n,v}$ for its cyclotomic
$\ZM_p$-extension. For all $n$ the natural map
$$(\NC_\infty)_{\Ga_n} \longrightarrow \NC_n$$ fits into an exact
sequence
$$\xymatrix{0 \ar[r] & (\NC_\infty)_{\Ga_n} \ar[r] & \NC_n\ar[rr]^-{\oplus_{v\in S_n} \text{\rm Artin at } v} &\quad
 & \bigoplus_{v\in S_n} \Gal((K_{n,v})^{\rm cyc}_\infty/K_{n,v})\ar[r]& 0.}$$
\end{prop}
\dem Let us examine the places above $p$ along $K_\infty/K$. These
primes have non trivial conjugated (hence equal) decomposition
subgroups. There exists a $d\in \NM$ such that $(i)$ no primes above
$p$ splits anymore in $K_\infty/K_{d}$ and $(ii)$ all primes
above $p$ are totally split in $K_{d}/K$. For all $n\geq d$ we
then have $S_n\cong S_{d}$, $\displaystyle \NC_n=\bigoplus_{v\in
S_{d}} L_{K_{n,v}}$ and $\Gal(K_\infty/K_n)$ identifies with the
local Galois groups $\Gal((K_{n,v})^{\rm cyc}_\infty/K_{n,v})$ (at
every $v\in S_{n}$). So in the case $n\geq d$, the proposition
follows from lemma \ref{locdes}. Next suppose that $0\leq n<d$,
and let $G^m_n$ be the Galois group $\Gal(K_m/K_n)$. Of course we
have $(\NC_\infty)_{\Ga_n}\cong((\NC_\infty)_{\Ga_{d}})_{G^{d}_n}$.
Using the previous case we have an exact sequence $$(\dag) \quad
\xymatrix{0 \ar[r] & (\NC_\infty)_{\Ga_{d}} \ar[r] &
\NC_{d}\ar[r]
 & \bigoplus_{v\in S_{d}} \Gal((K_{d,v})^{\rm cyc}_\infty/K_{d,v})\ar[r]& 0}$$
As $(K_{d,v})^{\rm cyc}_\infty/K_{d,v}$ is cyclotomic, the
action  of the local Galois group $\Gal(K_{d,v}/\QM_p)$ on the
group $\Gal((K_{d,v})^{\rm cyc}_\infty/K_{d,v})$ is trivial.
Hence we have an isomorphism of Galois modules
$$\bigoplus_{v\in S_{d}} \Gal((K_{d,v})^{\rm cyc}_\infty/K_{d,v})\cong \ZM_p[S_{d}].$$
Since the places above $p$  split totally in $K_{d}/K$, over  $G^{d}_0$ these
two modules are cohomologically trivial, and the same
is true for $\NC_{d}\cong \NC_0\otimes_{\ZM_p} \ZM_p[G^{d}_0]$
(as $G^{d}_0$-modules). Therefore in the sequence $(\dag)$ two
(hence all) terms are $G^{d}_0$-cohomologically trivial. The
triviality of $\widehat{H}^1(G^{d}_n,(\NC_\infty)_{\Ga_{d}})$
proves the injectivity of $(\NC_\infty)_{\Ga_n} \longrightarrow
\NC_n$. It then suffices to apply $N_{K_{d}/K_n}$ to the
sequence $(\dag)$ and use the triviality of the three
$\widehat{H}^0(G^{d}_n,--)$ to get the full sequence
$$\xymatrix{0 \ar[r] & (\NC_\infty)_{\Ga_n} \ar[r] &
\NC_n\ar[rr]^-{\oplus_{v\in S_n} \text{\rm Artin at } v} &\quad
 & \bigoplus_{v\in S_n} \Gal((K_{n,v})^{\rm cyc}_\infty/K_{n,v})\ar[r]& 0}.$$

\qed

Let $d$ be as before. Let $D_n$ be the decomposition subfield
for the place $p$ in $K_n$ (i.e $p$ is totally split in $D_n$ and
no place of $D_n$ above $p$  splits anymore in $K_n/K$). Then for
all $n\geq d$ we have $D_n=D_{d}$. Let us put $D=D_{d}$ for
the decomposition subfield for the place $p$ in $K_\infty$.
\begin{prop}\label{uslocdes} Recall that for a $\La$-module $M_\infty=\limpro(M_n)$ we denote
$\Ker_n (M_\infty)$ for the kernel, $\Mti_n$ for the image and
$\Coker_n (M_\infty)$ for the cokernel of the natural map $(M_\infty)_{\Ga_n}\longrightarrow M_n$.
\begin{enumerate}
\item For all $n\geq 0$, $\Ker_n(\UC_\infty)$ is isomorphic to $\ZM_p[S_n]$, hence its $\ZM_p$-rank is $\# S_n$.
\item For all $n\geq 0$, $\Coker_n(\UC_\infty)\cong N_{K_n/D_n}(\UC_n)$, hence its $\ZM_p$-rank is $\# S_n$.
In other words we  have an exact sequence $$\xymatrix {0\ar[r]
&\UCt_n \ar[r] & \UC_n \ar[r]^-{N_{K_n/D_n}} & \UC_{D_n} }.$$
\end{enumerate}
\end{prop}
\dem For any local field $E$ the (normalized) valuation $v$ of $E$
gives an exact sequence $0 \rightarrow U^1_{E} \rightarrow
\overline{E^\times}   \overset {v}  {\rightarrow} \ZM_p\rightarrow
0$, where $U^1_E$ is the set of principal units of $E$. Taking the
semi-local version of this sequence and projective limits we
obtain
$$\xymatrix{0\ar[r] & \UC_\infty \ar[r] & \NC_\infty \ar[r] & \ZM_p[S_{d}]\ar[r]& 0}.$$
But for all $n$, $\NC_n^{\Ga}=\NC_0$ contains no infinitely
$p$-divisible element. Therefore $\NC_\infty^\Ga=\{0\}$ and taking
$\Ga_n$-cohomology on this sequence we get :
$$\xymatrix{0\ar[r] & \ZM_p[S_{d}]^{\Ga_n} \ar[r] &
(\UC_\infty)_{\Ga_n} \ar[r] & (\NC_\infty)_{\Ga_n} \ar[r] & (\ZM_p[S_{d}])_{\Ga_n}\ar[r]& 0}
$$
Now since $(\NC_\infty)_{\Ga_n}\longrightarrow \NC_n$ is a monomorphism, $\Ker_n (\UC_\infty)$
identifies with $ \ZM_p[S_{d}]^{\Ga_n}\cong \ZM_p[S_n]$. This proves 1.

For 2. we use the notations $d$, $D_n$ and $D$ of the proof of
proposition \ref{slocdes}. We first prove the case $n\geq d$. By
compactness we have $\UCt_n=\bigcap_{m\geq n} N_{K_m/K_n}
(\UC_m)$. Because $n\geq d$ the global norm $N_{K_m/K_n}\colon
\UC_m\longrightarrow \UC_n$ is nothing else but the direct sum
place by place of the local norms of the extensions
$K_{m,w}/K_{n,v}$ (at each unique $w\in S_m$ above a fixed $v\in
S_n$). Fix a place $v\in S_n$ and for all $m\geq n$ still call $v$
the unique place of $S_m$ above $v$ and also $v$ the unique place
of $D_n$ under $v$. Note that $D_{n,v}=\QM_p$. Let $u\in
U^1_{K_{n,v}}$, then using local class field theory we see that
$u\in \bigcap_{m\geq n} N_{K_{m,v}/K_{n,v}} (U^1_{K_{m,v}})$ if
and only if $N_{K_{n,v}/\QM_p} (u)\in \bigcap_{n\in \NM}
N_{(\QM_p)^{\rm cyc}_n /\QM_p} (U^1_{(\QM_p)^{\rm cyc}_n})$, where
$(\QM_p)^{\rm cyc}_n$ denotes the $n^{\text{th}}$-step of the
cyclotomic (hence totally ramified) $\ZM_p$-extension of $\QM_p$.
By local class field theory $( ., (\QM_p)^{\rm cyc}_\infty/\QM_p)$ is an isomorphism from $U^1_{\QM_p}$ to $\Gal((\QM_p)^{\rm cyc}_\infty/\QM_p)$.
The equivalence
$$u\in \bigcap_{m\geq n}
N_{K_{m,v}/K_{n,v}} (U^1_{K_{m,v}}) \iff N_{K_{n,v}/D_{n,v} }
(u)=1$$ follows. As no place above $p$ splits in $K_n/D_n$ this gives
$$\UCt_n=\ker (N_{K_n/D_n} \colon \UC_n \longrightarrow
\UC_{D_n}).$$ Next pick an $n$ such that $0\leq n < d$. Consider
the diagram of fields
$$\xymatrix {K_{d} \ar@{-}[rd]\ar@{-}[d] & \\
K_n\ar@{-}[rd] & D_{d} \ar@{-}[d]\\ & D_n}$$ There every prime
above $p$ is totally split in the extensions $K_{d}/K_n$, and  $D_{d}/D_n$ and no prime
above $p$ splits at all in $K_n/D_n$ nor in $K_{d}/D_{d}$. It follows that
$N_{K_{d}/K_n}$ is surjective onto $\UC_n$ and that
$$\UCt_n=N_{K_{d}/K_n}(\UCt_{d})\subset \ker (N_{K_n/D_n}
\colon \UC_n \longrightarrow \UC_{D_n}).$$ Conversely pick
$u=(u_v)_{v\in S_n} \in \ker (N_{K_n/D_n} \colon \UC_n
\longrightarrow \UC_{D_n})$. At each $v\in S_n$ choose a single
$w(v)\in S_{d}$ above $v$ and define $t=(t_w)_{w\in S_{d}}\in
\UC_{d}$ by putting $t_w=1$ if there does
not exist $v$ such that $w= w(v)$ and $t_{w(v)}=u_v$ for all $v$ in $S_n$. Then we have $N_{K_{d}/K_n}
(t)=u$ and $t\in \UCt_{d}$ which shows that $u\in \UCt_n$.

\qed

Similar but not so precise statements about the sequence of global
units could be deduced from the proposition \ref{uslocdes} and
from the following proposition \ref{kuzmin}. As they are not
needed we don't state them. To end this section we recall the
analogous proposition for the sequence $\Uba^\pr_n$ of
$(p)$-units, which is a result of Kuz$ ^\pr$min (\cite{Ku72},
theorem 7.2 and theorem 7.3). Recall that $K$ is totally real, so
that $r_1=[K:\QM]$, $r_2=0$, and all $\U^\pr_n$ are
$\ZM_p$-torsion free. To avoid ugly notations $\widetilde {\Mba}$
we will denote $\Uti^\pr_n$ for the image of $\U^\pr_\infty$ in
$\U^\pr_n$.
\begin{prop}\label{kuzmin} \

\begin{enumerate}
\item The $\La$-module $\Uba_\infty^\pr$ is free of rank $[K:\QM]$.
\item For all $n$ the natural map $(\U^\pr_\infty)_{\Ga_n} \longrightarrow \U^\pr_n$ is injective.
\item $\Uti^\pr_n\cong (\U^\pr_\infty)_{\Ga_n}$ is a free $\ZM_p[G_n]$-module of rank $[K:\QM]$
(hence is $\ZM_p$ free of rank $p^n [K:\QM]=[K_n:\QM]$).
\end{enumerate}
\end{prop}
\dem 1 is theorem 7.2 and 2 is theorem 7.3 of \cite{Ku72}. There
the number field $K$ is arbitrary and a considerable amount of
effort is made to avoid using Leopoldt's conjecture. Another
proof (also without using Leopoldt's conjecture) is in \cite{KNF}. On the
other hand our abelian number field $K$ does satisfy Leopoldt's
conjecture so we may give, for the convenience of the reader, the
following shorter proof. We may suppose $n=0$ (else replace $K$ by
$K_n$). Let $\XG_\infty$ be the standard Iwasawa module
$\XG_\infty=\limpro(\XG_n)$,  where $\XG_n$ is the Galois group
over $K_n$ of its maximal abelian $S$-ramified $p$-extension (Nota
: of course this definition works also for $n=\infty$). Recall
from the end of the introduction the notation $X^\pr_n$ for the
$p$-part of the $(p)$-class group of $K_n$. From class field
theory we have the (some time called decomposition) exact sequence
$$\xymatrix{0\ar[r] & \U^\pr_\infty \ar[r] & \NC_\infty \ar[r] &
\XG_\infty \ar[r]& X_\infty^\pr \ar[r] & 0 }.$$ Put $\DC_\infty
=\im(\NC_\infty\longrightarrow \XG_\infty)$. By Leopoldt's
conjecture for $K$ we have $\XG_\infty^\Ga=0$ therefore
$\DC_\infty^\Ga=0$. This implies that the induced map
$(\U^\pr_\infty)_\Ga\longrightarrow (\NC_\infty)_\Ga$ is a
monomorphism. Now 2. follows from proposition \ref{slocdes}.

$K$ is totally real and $p\neq 2$, so  $(\U^\pr_\infty)_\Ga$  is
$\ZM_p$-free as a submodule of $\U^\pr_0$. As for $\NC_\infty$, we
have $(\U^\pr_\infty)^\Ga=0$ (same argument applies). These two
facts suffice to show the $\La$-freeness of $\U^\pr_\infty$. To
compute the rank we consider the sequence
$$\xymatrix {0\ar[r] & (\U^\pr_\infty)_\Ga\ar[r] & (\NC_\infty)_\Ga \ar[r] & (\DC_\infty)_\Ga \ar[r]& 0}.$$
As $\DC_\infty$ is a torsion-$\La$-module (as a submodule of
$\XG_\infty$) with trivial $\Ga$-invariants, its
$\Ga$-coinvariants are finite. Hence $(\U^\pr_\infty)_\Ga$  has
the same $\ZM_p$-rank as $(\NC_\infty)_\Ga$. By proposition
\ref{slocdes} this rank is $\rm{rank}_{\ZM_p}(\NC_0)-\#
S_0=[K:\QM]$. By Nakayama's lemma the $\La$-rank of
$\U^\pr_\infty$ is also $[K:\QM]$. This concludes the proof of 1.

3 is an immediate corollary of 1 and 2.

\qed

\section{From semi-local to global and vice-versa}

We now state and proceed to prove our main result in this paper. We
want to show that descent works asymptotically well for
$M_\infty=\Uba_\infty/\C_\infty$ or (which will be proven to be
equivalent) for $M_\infty=\UC^{(0)}_\infty/\C_\infty$ (see
explanations and notations below for the symbol $^{(0)}$). This
means that in both cases above, $\ker_n(M_\infty)$ and
$\coker_n(M_\infty)$ are finite of bounded orders. Our strategy of
proof is the  following. First we will show that bounding the
kernels and cokernels associated to both modules is equivalent.
Then we will use the injectivity of descent on $\Uba^\pr_\infty$
(proposition \ref{kuzmin}) to bound the kernels of descent for
$\U_\infty/\C_\infty$. Then we use local class field theory
(proposition \ref{uslocdes}) to bound the cokernels of descent for
$\UC^{(0)}_\infty/\C_\infty$.

But before this we have to slightly change the sequence
$\UC_n/\C_n$. Indeed, for all $n$, the module $\UC_n/\C_n$ is (by
Leopoldt's conjecture) of $\ZM_p$-rank $1$ while
$(\UC_\infty/\C_\infty)_{\Ga_n}$ is torsion. This rank $1$ comes
from $N_{K_n/\QM}(\U_n)=\{0\}$, or if we adopt the class field
theory point of view, it represents the rank of the
$\ZM_p$-extension $K_\infty/K$. Let us put the

{\bf Notation:} Let $F$ be a number field. In the sequel
$\UC^{(0)}_F$ will denote the kernel of $$N_{F/\QM}\colon \UC_F
\longrightarrow \UC_\QM .$$ Consistently $\UC^{(0)}_n$ will denote
the kernel of $N_{K_n/\QM}$ and $\UC_\infty^{(0)}=\limpro
\UC_n^{(0)}$.

By proposition \ref{uslocdes} we have $\UCt_n\subset \UC^{(0)}_n$
and therefore $\UC_\infty=\limpro(\UC^{(0)}_n)$.  Moreover
$\UC^{(0)}_n/\C_n$ is torsion and since $N_{K_n/\QM} (\UC_n)$ is
$\ZM_p$-torsion free we have $\UC^{(0)}_n/\C_n=\tor_{\ZM_p}
(\UC_n/\C_n)$, and accordingly $\UC^{(0)}_n/\U_n=\tor_{\ZM_p}
(\UC_n/\U_n)$. For all these reasons, it is clearly more
convenient to (and we will from now)  use the sequence
$(\UC_n^{(0)})_{n\in\NM}$ instead of $(\UC_n)_{n\in\NM}$. This
convention gives sense to the notations $\UCt^{(0)}_n$,
$\ker_n(\UC_\infty^{(0)})$, $\coker_n(\UC_\infty^{(0)})$ and
consistently to the same notations associated to sequences
$\UC^{(0)}_n/\C_n$, $\UC^{(0)}_n/\U_n$ and so on. Of course, due
to proposition \ref{uslocdes}, only the various $\coker_n$ will
actually change when we replace $\UC_n$ by $\UC_n^{(0)}$.
\begin{prop}\label{dna} Recall that $X_n$ stands for the $p$-part of the class group of $K_n$, and that
$\ker_n (X_\infty)$ is the kernel of the natural map $(X_\infty)_{\Ga_n}\longrightarrow X_n$.
\begin{enumerate}
\item The map $(\UC_\infty^{(0)}/\U_\infty)_{\Ga_n} \longrightarrow \UC^{(0)}_n/\U_n$ fits into
an exact sequence $$\xymatrix{ 0\ar[r]& (X_\infty)^{\Ga_n}\ar[r] &
(\UC_\infty^{(0)}/\U_\infty)_{\Ga_n} \ar[r] & \UC^{(0)}_n/\U_n
\ar[r] & \ker_n (X_\infty) \ar[r] & 0}$$
\item $\ker_n(\UC_\infty^{(0)}/\U_\infty)$ and $\coker_n(\UC_\infty^{(0)}/\U_\infty)$ are finite
and of bounded orders.
\end{enumerate}
\end{prop}
\dem By global class field theory we have an (some time called
inertia) exact sequence~:
$${(\rm R)} \quad \quad \quad
\xymatrix{0\ar[r] & \UC_n/\U_n \ar[r] & \XG_n\ar[r] & X_n \ar[r] & 0.}$$
Taking the $\ZM_p$-torsion counterpart of the sequence (R)  we have :
$$\xymatrix{0\ar[r] & \UC_n^{(0)}/\U_n \ar[r] & \tor_{\ZM_p}(\XG_n)\ar[r] & X_n .}$$
On the other hand if we take limits up to $K_\infty$ on (R), then
apply $\Ga_n$-cohomology, we obtain (by Leopoldt,
$\XG_\infty^{\Ga_n}=0$ and therefore $X_\infty^{\Ga_n}$ is finite):
$$\xymatrix{0\ar[r] &X_\infty^{\Ga_n} \ar[r]
& (\UC_\infty^{(0)}/\U_\infty)_{\Ga_n} \ar[r] &
(\XG_\infty)_{\Ga_n} \ar[r] & (X_\infty)_{\Ga_n}\ar[r] & 0.}$$
Lemma \ref{locdes}  has a global analogue which is the following :
\begin{lem}\label{stades}
The natural map $(\XG_\infty)_{\Ga_n} \longrightarrow \XG_n$ fits into the exact sequence
$$\xymatrix{0\ar[r]& (\XG_\infty)_{\Ga_n} \ar[r]& \XG_n\ar[r]^-{\rm res}& \Gal(K_\infty /K_n)\ar[r]& 0.}$$
In other words, descent provides an isomorphism
$$(\XG_\infty)_{\Ga_n} \cong \tor_{\ZM_p} \XG_n .$$
\end{lem}
\dem The proof of $(\XG_\infty)_{\Ga_n} \hookrightarrow \XG_n$ is exactly
the same as for \ref{locdes} : we only need to replace the $\GC_n$
there by the group $G_S(K_n)$ which is the Galois group of the
maximal $S$-ramified extension of $K_n$. The remaining part of the
exact sequence comes from maximality properties defining
$(\XG_\infty)_{\Ga_n}$ and $\XG_n$. Let $M_n$ be the maximal
abelian $S$ ramified $p$-extension of $K_n$. Then, by Leopoldt's
conjecture, one has $\tor_{\ZM_p} (\XG_n)=\Gal(M_n/K_\infty)$,
which gives the isomorphism. \qed

Let us resume the proof of 1 of \ref{dna}. Putting together the
three last sequences we obtain
$$\xymatrix{ & &0\ar[d] & & \\ (X_\infty)^{\Ga_n} \ar@{^{(}->}[r]
 & (\UC_\infty^{(0)}/\U_\infty)_{\Ga_n} \ar[d] \ar[r] & (\XG_\infty)_{\Ga_n} \ar[r] \ar[d] &
 (X_\infty)_{\Ga_n}\ar[r]\ar[d] & 0 \\
 0\ar[r] & \UC_n^{(0)}/\U_n \ar[r] & \tor_{\ZM_p}(\XG_n)\ar[r]\ar[d] & X_n \\ & &0 &
 &}$$
1 of \ref{dna} follows then from the snake lemma.

By Leopoldt's conjecture the $(X_\infty)^{\Ga_n}$'s are finite. As
$X_\infty$ is a noetherian $\La$-module the ascending union
$\bigcup_{n\in\NM} (X_\infty)^{\Ga_n}$ stabilizes. This shows that
the orders of $\ker_n(\UC_\infty^{(0)}/\U_\infty)$ stabilize. As
for $\coker_n(\UC_\infty^{(0)}/\U_\infty)\cong\Ker_n(X_\infty)$,
the maps $X_{n+1}\longrightarrow X_n$ (and consequently
$X_\infty\longrightarrow X_n$) are surjective as soon as
$K_{n+1}/K_n$ does ramify. By the classical Iwasawa theorems (see
\S 13 of \cite{Wa}), the orders of $(X_n)$ is asymptotically
equivalent to $p^{\la_X n + \mu_X p^n}$, where $\la_X$ and $\mu_X$
are the structural invariants of $X_\infty$. Since the order of
$(X_\infty)^\Ga$ is finite the orders of $(X_\infty)_{\Ga_n}$ are
also asymptotically equivalent to $p^{\la_X n + \mu_X p^n}$. This
proves that the orders of $\Ker_n(X_\infty)$, hence those of
$\coker_n(\UC_\infty^{(0)}/\U_\infty)$, are bounded and concludes
the proof of 2.

\qed

By Leopoldt's Conjecture we have for all $n$ :
$(\UC^{(0)}_\infty/\U_\infty)^{\Ga_n}\subset\XG_\infty^{\Ga_n}=0$.
And therefore an exact sequence
$$(\U_\infty/\C_\infty)_{\Ga_n}\hookrightarrow
(\UC^{(0)}_\infty/\C_\infty)_{\Ga_n} \twoheadrightarrow
(\UC^{(0)}_\infty/\U_\infty)_{\Ga_n}.$$ By the snake lemma with
the analogue exact sequence at  finite level we obtain the
sequence $$\xymatrix{0\ar[r] & \ker_n ( \U_\infty/\C_\infty)
\ar[r] & \ker_n ( \UC^{(0)}_\infty/\C_\infty) \ar[r] &
\ker_n(\UC^{(0)}_\infty/\U_\infty) \ar `r[d] `d[ld] `[ldldl] `d[d]
[lldd]& \\ & & & & \\ & \coker_n(\U_\infty/\C_\infty)\ar[r] &
\coker_n(\UC^{(0)}_\infty/\C_\infty)\ar[r] &
\coker_n(\UC^{(0)}_\infty/\U_\infty)\ar[r]& 0.}$$ Now, using this
sequence and the proposition \ref{dna} we prove our first key
lemma :
\begin{lem}\label{key1}\

\begin{enumerate}
\item $\Ker_n(\U_\infty/\C_\infty)$ and $\Ker_n(\UC^{(0)}_\infty/\C_\infty)$ are finite.
Their orders are simultaneously bounded or not.
\item  $\coker_n(\U_\infty/\C_\infty)$ and $\coker_n(\UC^{(0)}_\infty/\C_\infty)$ are finite.
Their orders are simultaneously bounded or not.
\end{enumerate}
\qed
\end{lem}

\begin{itemize}
\item[] {\it Remark :}
We have proven that, even if not bounded, the sequences of orders
would have been asymptotically equivalent. We will not use this,
because we will now proceed in bounding these orders !
\end{itemize}

\section{Descent kernels}

The second key lemma is
\begin{lem}\label{kerdes}\
\begin{enumerate}
\item
The orders of the kernels of the natural maps
$ (\Uba_\infty/\Cba_\infty)_{\Ga_n}\longrightarrow \Uba_n/\Cba_n$
are bounded.
\item The orders of the kernels of the natural maps
$ (\UC^{(0)}_\infty/\Cba_\infty)_{\Ga_n}\longrightarrow \UC^{(0)}_n/\Cba_n$
are bounded.
\end{enumerate}
\end{lem}
\dem
 From the commutative diagram :
$$\xymatrix{ & (\Cba_\infty)_{\Ga_n}\ar[r]\ar[d] & (\Uba_\infty)_{\Ga_n}
\ar[r]\ar[d] & (\Uba_\infty/\Cba_\infty)_{\Ga_n} \ar[r]\ar[d] & 0 \\
0\ar[r] &\Cba_n\ar[r] & \Uba_n\ar[r] & \Uba_n/\Cba_n \ar[r] &0 \\}$$

we deduce the exact sequence :
$$\xymatrix{0\ar[r] &\frac { \ken (\Uba_\infty )}
{ \im(\ken(\Cba_\infty))}\ar[r]& \ken (\Uba_\infty/\Cba_\infty) \ar[r]
& \Cba_n/\Cti_n }$$
By the part 2 of theorem \ref{Leo} $\ken (\Uba_\infty/\Cba_\infty)$ is
finite.  To control $\Coker_n(\C_\infty)$ we use the lemma
\begin{lem}\label{dist}
There exists an $N$ such that for all $n\geq N$ we have
$\Cba_n/\Cti_n\cong \Cba_N/\Cti_N$
\end{lem}
\dem Let $I$ be the inertia subfield of $p$ for $K_\infty/\QM$. By
lemma 2.5 of \cite{pmb}, for $n$ large enough ($n$ such that
$I\subset K_n$ is large enough), we have $\C_n=\Cti_n C_I$. It
follows that $\Coker_n(\C_\infty)\cong C_I/(\Cti_n\bigcap C_I)$.
Now the increasing sequence $\Cti_n\bigcap \C_I$ has to stabilize
because $\C_I$ is of finite $\ZM_p$-rank. This shows the lemma
\ref{dist}

\qed

To prove the lemma \ref{kerdes} it then suffices to bound the
orders of $$\Ken (\Uba_\infty ) / \im(\Ken (\Cba_\infty)).$$
For that, we prove that this sequence of quotients stabilizes :
\begin{prop}
For any noetherian $\La$-module $M_\infty$ let $\MI (M_\infty)$
denote the submodule of $M_\infty$ defined as follows :
$$\MI (M_\infty)=\bigcup_{n\in\NM} (M_\infty)^{\Ga_n}.$$ Exists $N$ such that for all $n\geq N$ we have
$$\Ken (\Uba_\infty ) / \im(\Ken (\Cba_\infty))\cong \frac {\MI(\U^\pr_\infty/\U_\infty)}
{\im (\MI(\U^\pr_\infty/\C_\infty))}$$
\end{prop}
\dem Starting from the sequence :
$$ \xymatrix{0\ar[r] &\Cba_\infty \ar[r] & \Uba^\pr_\infty \ar[r] &
\Uba_\infty^\pr/\Cba_\infty \ar[r] & 0 }
$$
\noindent one gets the diagram
$$\xymatrix{ 0\ar[r]& (\Uba^\pr_\infty/\Cba_\infty)^{\Ga_n} \ar[r]\ar[d]
& (\Cba_\infty)_{\Ga_n} \ar[r] \ar[d] & (\Uba^\pr_\infty)_{\Ga_n}
 \ar[d] \\ & 0\ar[r] & \Cba_n \ar[r] & \Uba^\pr_n \\} $$

By the proposition \ref{kuzmin}  the kernels
$\ken(\Uba^\prime_\infty)$ are trivial. Hence we have an
isomorphism $\ken(\Cba_\infty) \simeq (\Uba^\pr_\infty
/\Cba_\infty)^{\Ga_n}$. Since $\U^\pr_\infty$ is a noetherian
module, so is $\Uba^\pr_\infty /\Cba_\infty$ and therefore the
increasing sequence $(\Uba^\pr_\infty /\Cba_\infty)^{\Ga_n}$
stabilizes and for $n$ large we have $\ken(\Cba_\infty)\simeq
\MI(\U^\pr_\infty/\C_\infty)$. The same arguments  with
$\Uba_\infty$ instead of $\Cba_\infty$ proves that (provided $n$
greater than some $N$) we have $\ken(\Uba_\infty)\cong
\MI(\U^\pr_\infty/\U_\infty)$. This shows the proposition and also
(putting everything together) the first part of lemma
\ref{kerdes}. The second part of \ref{kerdes} follows then from
lemma \ref{key1}.

\qed

\section{Descent cokernels}

The third and final key lemma is

\begin{lem}\label{cokerdes}\
\begin{enumerate}
\item
The orders of the cokernels of the natural maps
$ (\Uba_\infty/\Cba_\infty)_{\Ga_n}\longrightarrow \Uba_n/\Cba_n$
are bounded.
\item The orders of the cokernels of the natural maps
$ (\UC^{(0)}_\infty/\Cba_\infty)_{\Ga_n}\longrightarrow \UC^{(0)}_n/\Cba_n$
are bounded.
\end{enumerate}
\end{lem}
\dem
Starting with the snake diagram
$$\xymatrix{ & (\Cba_\infty)_{\Ga_n} \ar[r] \ar[d]& (\UC_\infty)_{\Ga_n} \ar[r]\ar[d] &
(\UC_\infty/\Cba_\infty)_{\Ga_n} \ar[r]\ar[d] & 0 \\
0\ar[r] & \Cba_n \ar[r] & \UC_n^{(0)} \ar[r] &
\UC_n^{(0)}/\Cba_n\ar[r] & 0 }$$ one gets the sequence
$$\xymatrix{\ker_n(\UC^{(0)}_\infty/\Cba_\infty)\ar[r]& \Cba_n/\widetilde{C}_n \ar[r] &
\UC_n^{(0)}/\widetilde{\UC}_n \ar[r] &
\coker_n(\UC^{(0)}_\infty/\Cba_\infty)\ar[r] & 0. }$$ Recall that
$D$ is the maximal subfield of $K_\infty$ such that $p$ is totally
split in $D$. We may assume without loss of generality that
$D\subset K_n$ (else enlarge $n$). By lemma \ref{key1},
$\ker_n(\UC^{(0)}_\infty/\Cba_\infty)$ is bounded. By proposition
\ref{uslocdes} we have  $\widetilde\UC_n= \ker
(N_{K_n/D}\colon\UC_n^{(0)}\longrightarrow \UC_D)$. One then gets
an isomorphism $\UC_n^{(0)}/\widetilde{\UC}_n\cong
N_{K_n/D}(\UC_n^{(0)})$ and using this isomorphism the preceding
sequence reads

$$\xymatrix{ 0 \ar[r]& N_{K_n/D}(\Cba_n) \ar[r] & N_{K_n/D}(\UC_n^{(0)})
\ar[r] & \coker_n(\UC^{(0)}_\infty/\Cba_\infty)\ar[r] & 0. }$$
Now, 2 in lemma \ref{cokerdes} follows from the :
\begin{lem}\label{lcft} Let $c_n$ denotes $[K_n:D]$ (asymptotically $c_n$ is equivalent to $p^n$).
\begin{enumerate}
\item  $(\Cba_D)^{c_n}$ is a submodule of bounded finite index in $N_{K_n/D}(\Cba_n)$
\item $(\UC^{(0)}_D)^{c_n}$ is a submodule of bounded finite index in $ N_{K_n/D}(\UC_n^{(0)})$
\item $ (\UC^{(0)}_D)^{c_n}/\Cba_D^{c_n}$ is asymptotically equivalent to
$ \UC^{(0)}_D/\Cba_D$.
\end{enumerate}
\end{lem}
\dem Assertion 3 is immediate. The finite constant group
$\UC^{(0)}_D/\Cba_D$ maps onto
$(\UC^{(0)}_D)^{c_n}/(\Cba_D)^{c_n}$. Since the norm $N_{K_n/D}$
acts as $c_n$ on $\Cba_D$ and $\UC^{(0)}_D$ themselves, the
inclusions and finiteness of indices in 1 and 2 are clear. We have
to show that these finite indices are bounded.

For assertion 1 we use again lemma 2.5 of \cite{pmb}, that is
$\Cba_n=\widetilde C_n \Cba_I$. Moreover, as $\Cti_n\subset
\UCt_n^{(0)}$ we have $N_{K_n/D}(\widetilde C_n)=0$ by proposition
\ref{uslocdes} (Without using semi-local units,
$N_{K_n/D}(\widetilde C_n)=0$ can be checked directly using
distribution relations on a generating system of $\Cti_n$). Hence
we get $N_{K_n/D} (\Cba_n)=N_{K_n/D}(\widetilde C_n
\Cba_I)=N_{K_n/D}(\Cba_I)= N_{I/D}(\Cba_I)^{[K_n:I]}$. This gives
1 because, as $c_n$ itself, $[K_n:I]$ is asymptotically equivalent
to $p^n$ and $N_{I/D}(\Cba_I)$ is of (constant) finite index in
$\Cba_D$.

Assertion 2 is an easy exercise using local class field theory. Indeed,
recall that $\UC_n=\oplus_{v\mid p} \Uba^{1}_{K_{n,v}}$. Then the
global norm $N_{K_n/D}$ acts on each summand as the local norm
$N_{K_{n,v}/\QM_p}$. By local class field theory, the quotient
$\Uba^1_{\QM_p}/N_{K_{n,v}/\QM_p}(\Uba^{1}_{K_{n,v}})$ is
isomorphic to the $p$-part of the ramification subgroup of
$\Gal(K_{n,v}/\QM_p)$. These wild ramification subgroups are
cyclic with orders asymptotically equivalent to $p^n$. Summing up,
it follows that $N_{K_n/D}(\UC_n)$ contains $\UC_D^{c_n}$ with
bounded finite index. A fortiori $N_{K_n/D}(\UC_n^{(0)})$ contains
$(\UC^{(0)}_D)^{c_n}$ with bounded finite index. This concludes
the proof of lemma \ref{lcft} and therefore of the second claim in
lemma \ref{cokerdes}. The first  claim in \ref{cokerdes} then
follows  from lemma \ref{key1}.

\qed

With lemmas \ref{kerdes} and \ref{cokerdes} we have fullfilled our goal.
We have directly proved that the natural descent homomorphisms $(\U_\infty/\C_\infty)_{\Ga_n}\longrightarrow \U_n/\C_n$ have bounded kernels and cokernels. As a consequence, we get {\it without using Iwasawa's Main Conjecture} nor Sinnott's Index Formula an analogue for unit classes of Iwasawa's theorem.
Recall that any torsion $\La$-module $M_\infty$  has an invariant $\la$ which is the Weierstra\ss\ degree of (any) generators of its characteristic ideal and an invariant $\mu$ which is the maximal power of $p$ dividing (any) generators of its characteristic ideal. By purely abstract algebra it is classical and easy to prove that {\it if they are finite} the orders of $(M_\infty)_{\Ga_n}$ are asymptotically equivalent to $p^{\la n+\mu p^n}$.
\begin{theo}\label{iwaclassunit}
\
\begin{enumerate}
\item Let $\la_1$  and $\mu_1$ denote the structural invariants of
the Iwasawa module $\Uba_\infty/\Cba_\infty$. Then the orders of $\Uba_n/\Cba_n$ are
asymptotically equivalent to $p^{\la_1 n+\mu_1 p^n}$.
\item Let $\la_2$ and $\mu_2$ denotes the structural invariants of
the Iwasawa module $\UC_\infty/\Cba_\infty$. Then the orders of  $\UC^{(0)}_n/\Cba_n$ are asymptotically equivalent to $p^{\la_2 n+\mu_2 p^n}$.
\end{enumerate}
\end{theo}
\dem By propositions \ref{kerdes} and \ref{cokerdes}
the orders of $\Uba_n/\Cba_n$ are asymptotically equivalent to the orders of
$(\Uba_\infty/\Cba_\infty)_{\Ga_n}$. As they are finite the last orders are
equivalent to what we need. This shows 1. Same argument proves 2 as well.

\qed

\section{Two applications of Iwasawa's theorem for unit classes}\label{applications}
Our first application is a structural link between unit and ideal classes at infinity.
\begin{theo}\label{samelambdamu}\

\begin{enumerate}
\item The $\La$-modules $\U_\infty/\Cba_\infty$ and $X_\infty$ share the same structural
invariants $\la_1=\la_X$ and $\mu_1=\mu_X$.
\item  The $\La$-modules $\UC_\infty^{(0)}/\Cba_\infty$ and $\XG_\infty$ share the same structural
invariants $\la_2=\la_\XG$ and $\mu_2=\mu_\XG$.
\end{enumerate}
\end{theo}

\dem Consider the exact sequence of
$\La$-torsion modules
$$({\rm DNA})\quad \quad \quad \xymatrix{0\ar[r] & \U_\infty/\C_\infty \ar[r] & \UC^{(0)}_\infty/\C_\infty \ar[r] &
\XG_\infty\ar[r] & X_\infty \ar[r] & 0.}$$
The invariants $\la$ and $\mu$ are additive in exact
sequences. Therefore, going through the above (DNA) sequence, we
see that 1 is equivalent to 2. Now, by Sinnott's index formula, the orders of $X_n$ are
asymptotically equivalent to the orders of $\U_n/\C_n$. Using theorem \ref{iwaclassunit} and Iwasawa's theorem we get that the orders of $(X_\infty)_{\Ga_n}$ and $(\U_\infty/\C_\infty)_{\Ga_n}$ are finite and asymptotically equivalent. Therefore the sequence $\la_1 n +\mu_1 p^n$ is equivalent to the sequence $\la_X n +\mu_X p^n$ : assertion 1 follows.

\qed

As explained in the introduction, theorem \ref{samelambdamu}
is an immediate consequence
of the Iwasawa Main Conjecture. However the point here is that we
achieved a direct proof by only making use of Iwasawa classical
theorem, Sinnott's index formula, Coleman's morphism, and
Leopoldt's conjecture (the last two via theorem \ref{Leo}).

Conversely, theorem \ref{samelambdamu} could be used to simplify the (now classical) proof of the Main
Conjecture via Euler systems and the "class number trick". Let us recall the main lines, skipping technical details on characters. The $p$-adic $L$-functions are related via Coleman's theory to the characteristic series of $\UC^{(0)}_\infty/\C_\infty$, and one version of the Main Conjecture asserts that
$\UC^{(0)}_\infty/\C_\infty$ and $\XG_\infty$ have the same characteristic series (up to power of $p$, and
$\th$-componentwise for all Dirichlet characters $\th$ of the first kind). This is done in two steps~:
\begin{itemize}
\item[-] Use the Euler System of cicular units to "bound class groups" and to
show that the characteristic series of $X_\infty$ divide that of $\U_\infty/\C_\infty$. Hence, following the sequence $(DNA)$, to show that the characteristic series of $\XG_\infty$ divides that of $\UC^{(0)}_\infty/\C_\infty$. For full details, see \cite{Grei92}.
\item[-] To show the converse property, it suffices to prove the equality of the relevant $\la$-invariants. For this, one uses the "class number trick". Suppose that $K$ is the maximal real subfield of $M:=K(\ze_p)$, which is not a loss of generality. By Kummer duality the Iwasawa invariants of $\XG_\infty(K)$ and of $X_\infty^-(M)$ are equals. Then Iwasawa's asymptotical formula and the (minus part) of the class number formula show what we want.
\end{itemize}
It is this "class number trick" which could be advantageously replaced by theorem \ref{samelambdamu}.

Up to now we did not use Kummer duality nor any knowledge about the minus part of class groups, nor Ferrero-Washington's Theorem.
We use them now to show the vanishing of the $\mu$-invariant and thus remove the assertions "up to power of $p$" in all the discussions just above.

\begin{theo}\label{mu}
\
\begin{enumerate}
\item The structural invariant $\mu$ of the module
$\Uba_\infty/\Cba_\infty$ is trivial.
\item  The structural invariant $\mu$ of the module $\UC_\infty^{(0)}/\Cba_\infty$ is trivial.
\end{enumerate}
\end{theo}
\dem
%
We claim that all four modules in the above (DNA) sequence
have trivial $\mu$-invariant and we only need to
prove it for three out of them (actually two well chosen would be
enough). By 1 of theorem \ref{samelambdamu} $\mu_1=\mu_X$ and by theorem 7.15 of \cite{Wa} $\mu_X=0$.
Let us draw the main lines of the proof written in \cite{Wa} of the triviality of $\mu_X$.
\begin{itemize}
\item[step 1] The main ingredient is  Ferrero-Washington's theorem  \cite{FW} which claims that the power series
$g_\psi(T)$ of our first section is prime to $p$.
\item[step 2] Over $K(\ze_p)$, using step 1 and the analytic class number formula for the minus part, one
deduces that the sequence of orders of $X_n^-$ is equivalent to
$p^{\la n}$, where $\la$ is the Weierstrass degree of the product
of relevant $g_\psi(T)$'s. By Iwasawa's theorem, this implies the
triviality of the structural $\mu$-invariant of $X_\infty^-$.
\item[step 3] Using the classical mirror inequality $\mu^+\leq \mu^-$ derived from Kummer duality, one
recovers the triviality of the $\mu$-invariants of the plus part
$X^+_\infty$ over $K(\ze_p)$, which in turn implies the triviality
of the $\mu$-invariant of $X_\infty$ for our base field $K$.
\end{itemize}
For the remaining module $\XG_\infty$, a possible proof follows
the above first two steps by just replacing the analytic class
number formula for the minus part by Leopoldt's formula for the
order of the even part of the $\XG_n$ in terms of products of
values at $1$ of $p$-adic $L$-functions. Alternatively let us just
examine more carefully the third step. Actually Kummer duality (we
are only making use here of corollary 11.4.4 of \cite{NSW} but
full original Kummer duality is in \cite{I73}) gives that
$\mu(\XG_\infty^+)=\mu(X_\infty^-)$ over $K(\ze_p)$ and the
inequality in the third step follows from $\XG_\infty
\twoheadrightarrow X_\infty$. Hence step 2 gives directly
$\mu(\XG_\infty^+)=0$ over $K(\ze_p)$, which is all we need.
The triviality of $\mu_2$ follows and this concludes the proof of 2 of theorem \ref{mu}.
\par
\qed
\par
\begin{itemize}
\item[] {\it Remark :}
Assertion 1 of theorem \ref{mu} is also proven by Greither in the appendix of
\cite{FG04}. Greither's proof is slightly different because it
makes no use of Leopoldt's conjecture, and for that reason needs
to work with maybe  infinite $ |(\Uba_\infty/\Cba_\infty)_{\Ga_n}
|$. To deal with that difficulty, Greither introduces the notion
of "tame" sequences of modules (roughly speaking these are
sequences of modules whose inverse limits are without $\mu$).
\end{itemize}


\begin{thebibliography}{NSW}
\bibitem[Bt]{Tan} T. Beliaeva, \emph{Unit\'{e}s semi-locales modulo sommes de
Gau\ss\ en th\'eorie d'{I}wasawa}, Th\`{e}se de l'universit\'{e} de
Franche-Comt\'{e} Besan\c{c}on (2004).

\bibitem[Bjr]{pmb}
J.-R. Belliard, \emph{Sous-modules d'unit\'es en th\'eorie d'{I}wasawa},
  Th\'eorie des nombres, Ann\'ees 1998/2001, Publ. Math. UFR Sci. Tech. Besan\c
  con, 2002, 12~p.

\bibitem[FG]{FG04}
M. Flach, \emph{The equivariant {T}amagawa number conjecture: a survey},
  Stark's conjectures: recent work and new directions, Contemp. Math., vol.
  358, Amer. Math. Soc., Providence, RI, 2004, With an appendix by C.\
  Greither, pp.~79--125.

\bibitem[FW]{FW}
B.~Ferrero and L.~Washington, \emph{The {I}wasawa invariant $\mu_p$ vanishes
  for abelian number fields}, Ann. of Math. \textbf{109} (1979), 377--395.

\bibitem[G]{Grei92}
C. Greither, \emph{Class groups of abelian fields, and the main
  conjecture}, Ann. Inst. Fourier (Grenoble) \textbf{42} (1992), no.~3,
  449--499.

\bibitem[I1]{I64}
K. Iwasawa, \emph{On some modules in the theory of cyclotomic fields}, J.
  Math. Soc. Japan \textbf{16} (1964), 42--82.

\bibitem[I2]{I73}
K. Iwasawa, \emph{On $\mathbb{Z}_{\ell}$-extensions of algebraic
number
  fields}, Ann. of Math. (2) \textbf{98} (1973), 246--326.

\bibitem[KNF]{KNF}
M. Kolster, T.
Nguy{\cftil{e}}n Quang {\Dbar}{\cftil{o}}, and
  V. Fleckinger, \emph{Twisted ${S}$-units, $p$-adic class number
  formulas, and the {L}ichtenbaum conjectures}, Duke Math. J. \textbf{84}
  (1996), no.~3, 679--717.

\bibitem[K]{Ku72}
L.~V. Kuz$'$min, \emph{The {T}ate module of algebraic number fields}, Izv.
  Akad. Nauk SSSR Ser. Mat. \textbf{36} (1972), 267--327.

\bibitem[N1]{Ng84}
T. Nguy{\cftil{e}}n Quang {\Dbar}{\cftil{o}}, \emph{Formations de
  classes et modules d'{I}wasawa}, Number theory, Noordwijkerhout 1983
  (Noordwijkerhout, 1983), Springer, Berlin, 1984, pp.~167--185.

\bibitem[N2]{Ng05}
T. Nguy{\cftil{e}}n-Quang-{\Dbar}{\cftil{o}}, \emph{Sur la conjecture faible de Greenberg dans le cas ab\'elien $p$-d\'ecompos\'e},  Int. J. Number Theory \textbf{2} (1972), no.~1, 49--64.

\bibitem[NSW]{NSW}
J. Neukirch, A. Schmidt, and K. Wingberg, \emph{Cohomology of
  number fields}, Grundlehren der Mathematischen Wissenschaften, vol. 323, Springer-Verlag, Berlin,
  2000.

\bibitem[Si]{Si80}
W. Sinnott, \emph{On the {S}tickelberger ideal and the circular units of an
  abelian field}, Invent. Math. \textbf{62} (1980), no.~2, 181--234.

\bibitem[T]{T99}
T. Tsuji, \emph{Semi-local units modulo cyclotomic units}, J. Number Theory
  \textbf{78} (1999), no.~1, 1--26.

\bibitem[W]{Wa}
L.~C. Washington, \emph{Introduction to cyclotomic fields}, second ed.,
  Springer-Verlag, New York, 1997.
\end{thebibliography}

{\bf Acknowledgement :} I thank professor Th{\cfac{o}}ng Nguy{\cftil{e}}n Quang {\Dbar}{\cftil{o}}
for many instructive conversations and helpful comments on this paper.

\noindent Jean-Robert Belliard

\noindent belliard@math.univ-fcomte.fr

\noindent Laboratoire de Math\'ematiques de Besan\c{c}on

\noindent Universit\'e de Franche-Comt\'e

\noindent 16 Route de Gray

\noindent 25030 BESAN\c{C}ON Cedex
\end{document}